\documentclass[english]{amsart}

\usepackage{babel}

\usepackage{hyperref}
\hypersetup{
    colorlinks=true,
    linkcolor=red,
    filecolor=magenta,      
    urlcolor=cyan,
    citecolor=blue,
    pdftitle={Topomodal nominalism},
    }
\usepackage{mathptmx}
\usepackage{stmaryrd}
\usepackage{amssymb}
\usepackage{cleveref}

\theoremstyle{definition}

\newtheorem{axiom}{Axiom}
\crefname{axiom}{Axiom}{Axioms}

\crefname{remark}{Remark}{Remarks}
\newtheorem{theorem}{Theorem}

\title{An approach to topology via lattices of names}
\author{Martín Torres-Valverde}
\address{Department of Mathematics, Universitat Jaume I-IMAC, E-12071, Castell\'{o}, Spain}
\email{valverdm@uji.es}
\thanks{Research partially supported by the project CIACIF/2026/88 funded by Conselleria de Educaci\'o, Cultura, Universitats i Ocupaci\'o-Generalitat Valenciana, and by predoctoral program of the Universitat Jaume I PREDOC/2025/80.}



\begin{document}

\begin{abstract}
Synthetic topology encompasses various alternatives to the classical definition of a topological space by Hausdorff. In particular, Abstract Stone Duality (ASD) develops the theory of locally compact spaces directly from logic. This study introduces a semi-formal theory of synthetic topology using a first-order version of the intuitionistic modal logic iS4, inspired by a nominalist, dynamic interpretation of the specialization order. The remarkable similarities to ASD include a constructive interpretation of the open-closed duality, overtness of the line, compactness of the unit interval and the nuclei that define closed and open subspaces.
\end{abstract}

\maketitle

\textbf{Important note.} There were mistakes in the first version of the preprint, caused by careless use of dualities. As is seen here, the duality between closed and open subspaces is written as \(\tau(X) \cong \mathfrak f(X^*)^*\), whereas in the earlier pre-version it was written \(\tau(X) \cong \kappa(X^*)^*\), which is wrong. The solution was obtained by replacing the specialization order between spaces by a category, with \emph{parent-child relations} as morphisms, which are basically the same as the join-respecting functions \(\kappa(X) \to \kappa(Y)\) that were sometimes used in the first version. This is better, among other reasons, because it permits the definition of quotients.

This version corrects the mistakes above, but it shall take some time to recover some results. In particular, some of the connections to ASD still have to be worked out. Hence, one should expect an update in the following month.

\section{Introduction}

For many decades, it has been known that point-set topology is just one approach to topology among many, including variants of locale theory \cite{pointless}, formal topology \cite{formal-topology}, synthetic topology \cite{escardo-synthetic-topology}, and so on. Of particular relevance to this paper is Abstract Stone Duality (ASD), the system in \cite{TaylorP:dedras,TaylorP:lamcra}, that behaves in certain ways like classical systems (i.e., the open-closed duality) while retaining computational power. ASD models well the category of (computably based) locally compact spaces. However, the logic of ASD was intentionally chosen to be weak, as it was apparent that a strong logic would lose the strongly constructive interpretation. This makes it hard to reason informally in ASD.

This paper provides an expressive semi-formal theory, called TN,\footnote{TN is shorthand for \emph{Topomodal Nominalism}, a name based on its other inspiration: the \emph{modal nominalism} of \cite{weaver2025cc}.} that is closely linked to ASD, in the hope that it will make it easier to reason informally about topology independently of set theory. The basis of TN is a well-known idea in the area of semantics of computer languages: the \emph{specialization order} \cite{10.1007/BFb0012801,lipparini}. While the specialization order may be defined in an arbitrary topological space, it is discrete for classical Hausdorff spaces, which explains why its popularity has remained modest in areas such as analysis.

The specialization order can be alternatively regarded as an order between names. For example, the name ``A planet'' is more general than the name ``The planet Earth''; or equivalently, the name ``The planet Earth'' is more specific than the name ``A planet''. We see that a more general name is a name that refers to more objects. Each variable of TN stands for a name, and the specialization order between names is the most basic primitive. However, what a name refers to may change over time. This is an important fact in TN: a \emph{constant} is a name that, once bound, will always refer to the same object, and a \emph{filter} is a bound name that refers to increasingly specific objects. Importantly, there are only finitely many constants at each time, and even if there may be, in a sense, infinitely many variables, only finitely many of them are bound. The first part of the paper is dedicated to adapting the nominalist philosophy in \cite{weaver2025cc} to our purposes.

The other ingredient we need are \emph{spaces}. Spaces are variable, but finite, sets of objects. When is a space more specific than another? Several answers are valid. When more specific spaces are larger spaces, the resulting filters are \emph{open spaces}. By default, an object does \emph{not} belong to an open space, but it may end up there as more information is gained. On the other hand, when more specific spaces are smaller spaces, the resulting filters are called \emph{closed spaces}. By default, an object belongs to a closed space, but it may end up being outside the closed space after more information is gained. Hence, open and closed spaces are complementary.

In fact, there are multiple ways that a space can be more specific than another, and idea heavily inspired by \emph{Homotopy type theory} \cite{hottbook}, where equality can have multiple witnesses. The relation between TN and such systems should be worth exploring. Hence, equivalence of spaces is often taken to depend only on the structure of the spaces (their specialization order and some extra structure) and not on their particular elements. A choice of specialization order between spaces hence results in a category, called a \emph{type}.

Rather than the linear order of time, the specialization order evolves over a more general pre-order \((W, \le)\), which is is called the \emph{accessibility relation} between \emph{possible worlds}. We also have a ``holds'' relation (often called \emph{satisfaction})  \(w \Vdash \phi\) between worlds \(w \in W\) and sentences \(\phi\) such that \(w \Vdash \phi\) if and only if \(\phi\) holds at \(w.\) It is convenient to add two modalities `\(\square\)' \emph{necessarily} and `\(\Diamond\)' \emph{possibly} such that \(w \Vdash \square \phi\) if and only if, for all \(v \ge w,\) \(w \Vdash \phi,\) and \(\Diamond \phi\) if and only if there exists \(v \ge w\) such that \(v \Vdash \phi.\) The modalities of TN follow the rules of iS4, an intuitionistic modal logic that has been noted to accurately model topological questions: \cite{10.1093/logcom/exae030} show that iS4 is the ``logic of topological spaces, with persistence relative to the specialization order''.

The spaces that are expected to be well modelled in TN are locally compact spaces, given the similarities to ASD. In fact, each basic open in TN can dually be regarded as a \emph{part}, which has some of the properties of a compact subspace in classical topology. The sense in which TN models locally compact spaces has not been made precise yet, since the theory has not been turned into a formal system with minimal axioms.

\section*{Outline}
\begin{itemize}
    \item \Cref{sec:nominalism} is dedicated to reviewing the philosophy laid out in \cite{weaver2025cc}. It is seen how the modalities of iS4, the specialization order and judgements arise naturally when one takes the nominalist and constructivist philosophies seriously.
    \item \Cref{sec:specialization} expands on the specialization order.
    \item \Cref{sec:preliminaries} contains some preliminaries on finite order theory.
    \item \Cref{sec:logic} introduces the first-order version of the logic iS4 that TN uses.
    \item \Cref{sec:types} introduces types of various flavours.
    \item \Cref{sec:sierpinski} introduces the Sierpinski space.
    \item \Cref{sec:functions} introduces functions
    \item \Cref{sec:rationals} introduces the rational numbers.
    \item \Cref{sec:ascending} introduces the ascending and descending numbers.
    \item \Cref{sec:line} constructs the real line and extends addition of rational numbers to real numbers.
    \item \Cref{sec:applications} provides a glimpse on how locally compact mathematics could be developed in TN.
\end{itemize}

\section{Nominalism}
\label{sec:nominalism}

What is the nature of abstract objects such as numbers, sets or functions? This question has sparked debate for centuries, and there is no doubt that it will continue to do so. Rather than examining the 2 main competing philosophies (Platonism, nominalism) and their many variants, here we shall focus only on \emph{nominalism}, which states that abstract objects have no independent existence, and which is sometimes not taken seriously as a philosophy that can successfully explain mathematical practice.

Recently, Weaver's book \cite[Part II]{weaver2025cc} has proved that countabilist nominalist foundations are enough for developing the bulk of mathematics, while eliminating what can (perhaps controversially) be considered to be \emph{set-theoretic pathology}, such as the Stone-Čech compactification of the natural numbers or the duals of some non-separable Banach spaces. While weaker alternatives to set theory are not new, \cite[Part I]{weaver2025cc} includes novel treatments of the classical paradoxes based on the philosophy of \emph{truth constructivism}. However, the resulting theory \emph{Constructive Countabilism} makes compromises in order to recover most of mainstream mathematics. Here we want to describe \emph{new approaches} in even weaker, but still expressive, foundations, by unveiling the hidden modalities and specialization order in \cite[Part I]{weaver2025cc}.

\subsection{Mereological sums}
To begin with, \cite[ch. 2]{weaver2025cc} gives an answer to the question: what is a set? Starting with concrete objects (macroscopic, sharply defined physical objects) as basic, Weaver begins by distinguishing sets from \emph{mereological sums}:
\begin{quote}
    Let us begin by noting that one basic way to identify new concrete objects is as amalgamations of previously available concrete objects. The philosophical term for such a thing is ``mereological sum''. For instance, a deck of cards is the mereological sum of the individual cards, and each card is the mereological sum of the atoms which compromise it. \cite[p. 20]{weaver2025cc}
\end{quote}
Let's fix some notation for mereological sums. Given a name \(x\) and a name \(y,\) each of which stands for a physical object, we can form a name \(x \sqcup y\) (read `\(x\) join \(y\)') that refers to the mereological sum of \(x\) and \(y.\) For example, if \(x =\) ``the left brain hemisphere'' and \(y =\) ``the right brain hemisphere'' then \(x \sqcup y =\) ``the brain''. Or if \(x = \) ``apple 1'' and \(y = \) ``apple 2'' then \(x \sqcup y = \) ``apples 1 and 2''. We can also form joins of many names to refer to a mereological sum of many objects:
\[ \text{``deck''} = \bigsqcup_{x \text{ card}} x. \]

\subsection{Identity criteria}

Weaver continues by questioning what it means for two objects to be \emph{the same}:
\begin{quote}
    Is it really clear that the same object cannot simultaneously exist in distinct locations? If \emph{Gravity's Rainbow} may be found in fine bookstores everywhere, evidently it \emph{can} be in more than one place at once. The issue here is what we mean by ``the same'', what it means to say that those different copies of \emph{Gravity's Rainbow} are all the same book. Are you the same person you were ten years ago? Is a caterpillar that spins a cocoon the same creature as the butterfly that emerges later? \cite[p. 21]{weaver2025cc}
\end{quote}
He solves this question by introducing \emph{identity criteria}. A caterpillar is the same as the butterfly that emerges if we adopt the criterion that they are the same, and they are not the same if we adopt the criterion that they are not the same. Theseus' paradox is not a problem from this point of view:
\begin{quote}
    The ship of Theseus problem --- at what point is it no longer the same ship? --- is at heart a question about what identity criterion to use. Indeed, it is hardly a problem at all once we recognize that we are free to choose whatever identity criterion we like. Or rather, it becomes merely a problem about which identity criterion is most appropriate or useful in that setting. \cite[p. 22]{weaver2025cc}
\end{quote}
Here we shall say two names are \emph{equivalent}, \(x \equiv y,\) if both \(x\) and \(y\) refer to ``the same'' object. We can express the fact that ``identity criteria can change'' by saying that equivalence of names is \emph{contingent} on the identity criterion that is being used. Hence, even if \(x \equiv y,\) it may well be possible that \(x \not\equiv y.\) This can be expressed using the modality `\(\Diamond\)' (\emph{possibly}). If \(x = \) ``the butterfly'' and \(y = \) ``the caterpillar'', then we have:
\[ \Diamond (x \equiv y) \land \Diamond (x \not\equiv y). \]
Weaver himself chooses a different formalization of identity criteria, by having the variables \(x, y, z\dots\) range over the objects themselves rather than names. However, I like to add an extra layer of indirection and let the variables range over names of objects. This will also allow us to have names for objects that are not macroscopic, sharply-defined physical objects. It is safe to conflate an object with a \emph{constant} name that refers to it.

\subsection{Judgements}

Adopting a particular identity criterion is how a \emph{class} (in Weaver's terminology) or a \emph{space} (in our terminology) is constructed. The construction of a space is an example of \emph{judgement} \cite{martin1996meanings}. A judgement, for us, is simply an act that can be represented via the sentence ``this judgement has been made''. For example, in type theory, a variable \(x\) having the type \(X\) is a judgement, which we write as \(x : X\). As long as \(x\) is unbound, we can make the judgement \(x : X\) at any time, but the sentence ``\(x : X\)'' \emph{only} becomes true once we have made the corresponding judgement. Prior to the judgement, the sentence ``\(x : X\)'' was \emph{false}. The act of judgement takes place independently of the language we use to describe it, but we usually conflate a judgement with its linguistic representation.

Here we prefer to call the act itself a \emph{basic judgement}. A boolean combination of the corresponding sentences of type ``the basic judgement \(J\) has been made'' is called a \emph{judgement}. In this way, the conjunction, disjunction or negation of judgements is a judgement. However, we convene that, if a basic judgement has been made, then it \emph{necessarily} has been made; hence we may say basic judgements are \emph{open}. A world is \emph{accessible} (see \Cref{sec:logic}) if and only if more basic judgements have been made in it. While basic judgements may be undone, they can't be undone as we move up the accessibility relation, by definition.

If we take a variable `\(x\)' to denote a real number in the middle of a certain argument, then we will continue treating the variable \(x\) to be a real number in the subsequent steps; this is an example of an open judgement. However, let's adopt the following identity criterion: two judgements are the same if both have been made or neither has. (This is a version of the \emph{Sierpinski space}; in \Cref{sec:sierpinski}, we shall construct a simpler version.) This is a decidable identity criterion, and yet, in the course of an argument, two basic judgements may begin being `the same' (when neither has been made) only for them to become `different' (when one of them is made) only for them to become `the same again' (when the second one is made). Hence, this identity criterion is neither open nor closed; nevertheless, it is a judgement, as it is a boolean combination of basic judgements.

\subsection{Constructive truth as a judgement}

In \cite[ch. 4 -- ch. 6]{weaver2025cc}, Weaver introduces a variety of concepts including constructive truth (or \emph{assertibility}), meaningfulness, classical truth, concepts, classes, and so on. He advocates for the position that truth is constructed in stages, and provides a resolution for several classical paradoxes including the Liar paradox.

This is a well-thought take on truth; the only thing I can comment is that the modalities `\(\Diamond\)', `\(\square\)' appear implicitly, as early as in the choice of the word ``assertibility''. To assert something is a judgement, and something is assertible if it is possible to assert it. (Note that we are using the word ``possible'' in a strong sense in this paper: it is \emph{possible} to make a judgement if we are compelled to make it, not merely if we can conceivably make it.) Hence, modal logic provides a more fine-grained view of what assertibility is. Similarly, it is a judgement to consider a sentence meaningful.

\subsection{Finitism versus infinitism}

Weaver talks about how in mathematics one is concerned with ``things that are true in any world'' \cite[p. 27]{weaver2025cc}.\footnote{The modal language in the \emph{modal nominalism} described explicitly by Weaver is closer to S5, whereas our modal language is closer to S4 (or rather, its intuitionistic version iS4). In fact, the modalities are used for different purposes.} Then he distinguishes three views on the possibility of actual infinity in a possible world; here we will group countabilism and Cantorianism (as well as intermediate philosophies such as \emph{sequentialism}) under the name \emph{infinitism} \cite[p. 29]{weaver2025cc}:
\begin{itemize}
    \item According to \emph{finitism}, infinity is by definition potential: all sets in all possible worlds are finite.
    \item According to \emph{countabilism}, the position he advocates, there are some possible worlds with actually infinite sets, but none with uncountable sets. In particular, there is a canonical model of the natural numbers.
    \item According to \emph{Cantorianism}, there are possible worlds with uncountable sets; in fact, for every set actually existing in some world, there is some other world that contains its power-set.
\end{itemize}

Then Weaver introduces concrete \emph{proxies} for various abstract objects; in particular, a list of bound names is proxy for a set, with the identity criterion that two lists are considered \emph{the same} if they contain the names of the same objects (even if their elements vary in order and multiplicity) \cite[p. 34]{weaver2025cc}. This is a valid nominalist definition of a set, and perhaps it should be standard, as is his notion of a \emph{pure set} \cite[p. 37]{weaver2025cc}.

Weaver does not argue for countabilism until \cite[ch. 7]{weaver2025cc}. He acknowledges there isn't a ``knock-down argument for countabilism'', but that nevertheless it is a reasonable middle ground between finitism and Cantorianism \cite[p. 111]{weaver2025cc}. His argument basically is that actual infinite totalities (and infinitely-long computations) are conceivable; therefore infinitism should be sound and finitism too restrictive. I will not argue against the possibility of actual infinite sets.\footnote{I'm slightly sceptical that for \(\mathbb N\) to ``be in bi-univocal correspondence'' with a proper part of itself is the \emph{only} conceivable way that it can be actually infinite. To me, it appears just as likely that infinity is necessarily a \emph{relative} notion, so that what is infinite for us would be finite to an inhabitant of a ``locally infinite'' universe.}

My objection against countabilism and in favour of finitism (here understood broadly, encompassing the intuitionistic and other constructive schools) diverges a bit from the modal nominalism of Weaver, and it doesn't come from some extreme form of scepticism either. In my view, even if actually infinite sets did exist in some possible world (whatever that might mean), that does not mean that actual infinite sets should have a leading role in mathematical practice. Mathematics should concern itself \emph{primarily} with modelling \emph{our} world, and our world is, for most purposes, locally finite.

Note that topology is typically concerned with systems about which we have a finite amount of information. For example, of a real number \(x\) obtained through experimentation we can obtain some open information such as \(0.2 < x < 0.3,\) but we will never know all the digits of \(x\) (if that even made any sense). Were mathematics primarily concerned with studying the properties that are true in any world, we would not be especially interested in this topology for the real line. In fact, I expect that we would be more interested in \emph{some other} analogue space that inhabitants of a locally infinite universe would study empirically; something like a completion of the surreal numbers.

While infinitist systems have usually proved easier to work with, which makes the creation of infinitist systems that better explain what mathematicians are actually doing a worthy task, there remains unexplored terrain in finitist systems too. For example, the open-closed duality has been found in \emph{Abstract Stone Duality} (ASD) \cite{TaylorP:lamcra, TaylorP:dedras} not to be inherently unconstructive and not to require any pre-existing set theory (not even an intuitionistic one, even if models of ASD contain \emph{arithmetic universes}, which are models of a constructive theory of countable sets \cite{TaylorP:insema}). Hence, the rest of the paper is devoted to developing a theory TN where topology (hopefully, at least, the theory of locally compact spaces) can be developed constructively without some of the usual inconveniences of intuitionistic logic.

\section{The specialization order}
\label{sec:specialization}

While non-Hausdorff topological spaces are practically absent from the standard curriculum of a Mathematics degree, they happen to lie right at the intersection of mathematics and linguistics. See, for example, \cite{lipparini}, who creates an algebrization of topology via the specialization order related to TN (for example, morphisms are covariant with respect to the specialization order). When two points are not separated by opens, we may find that one of them is \emph{more special} or \emph{more specific} than the other (the more special one is contained in \emph{more} opens).\footnote{The multiple conventions for the specialization order are not reviewed here.} This \emph{specialization order} occurs no less than in language itself: a name \(x\) is said to be \emph{more special} than another name \(y,\) or \(y\) \emph{more general} than \(x,\) formally \(x \sqsubseteq y,\) if every object that \(x\) refers to is also an object that \(y\) refers to. For example:
\begin{quote}
    ``Screwdriver'' \(\sqsubseteq\) ``Tool'' \(\sqsubseteq\) ``Object''.
\end{quote}
When we have \(x \sqsubseteq y\) and \(y \sqsubseteq x,\) we say \(x\) and \(y\) are \emph{equivalent}, formally \(x \equiv y,\) coinciding with the definition given earlier. The specialization order is not new at all; however, its main applications have been developed in the study of \emph{computer languages} as opposed to natural languages, or even mathematical languages. In both classical and synthetic topology, an important application of these ideas is the \emph{Sierpinski space}, a space with two names, one more specific than the other, which classifies open (or closed) subspaces.

The specialization order is, in particular, connected to \emph{mereological sums}: if \(x \sqcup y\) is the join of the names \(x\) and \(y\) (which was defined to be a name that refers to the mereological sum of the objects that \(x\) and \(y\) refer to), then we have that \(x \sqcup y\) is actually the least upper bound of \(x\) and \(y\) in the specialization order. In other words, a name is more general than \(x \sqcup y\) if and only if it is more general than \(x\) and \(y.\)

We can also define the \emph{meet} \(x \sqcap y\) of two names to be their highest lower bound in the specialization order. For example,
\begin{quote}
    ``Something pink'' \(\sqcap\) ``A panther'' \(\equiv\) ``A pink panther''.
\end{quote}
We could define a name ``nothing'' (formally, \(\bot\)) that does not refer to any object, hence being more special than any other name, and a name ``something'' (formally, \(\top\)) that refers to any object, hence being the most general name. However, for our purposes, it makes more sense to take the specialization order as the only basic primitive. Minimal and maximal elements, as well as the lattice structure (joins and meets), will appear only in certain spaces.

How does all of this relate to topology? I actually believe there is \emph{little more} to topology than the specialization order of names. Even for a Hausdorff space, the specialization order appears in its lattice of opens. For a locally compact space in particular, its lattice of opens carries itself a natural topology called the \emph{Scott topology} \cite{escardo-exponential-spaces}, whose specialization order is given by (reverse) inclusion. If a space has a countable basis, then each basic open can be given a name.

We have shown that names form a (distributive) lattice under the specialization order. It is also clear that there can only be a finite amount of names bound at any time. Hence, the theory of finite distributive lattices will be the foundations of the theory TN presented here. However, language is not static: new names are defined, and the meaning of those that are already bound may vary over time. Whether a name refers to an object is a \emph{judgement}, because it is first of all an act and only secondarily a sentence. As we said, judgements are decidable, but we need the right theory to study the order in which judgements are formed and the way in which they vary over `time'. This is developed in the next sections.

\section{Order theory preliminaries}
\label{sec:preliminaries}

Let's review some results from finitist order theory. A \emph{space} is, for us, a finite collection \(X\) of \emph{variables}, equipped with an equivalence relation `\(\equiv\)' between variables, such that \(x \in X\) is decidable and \(x \equiv y\) is decidable. A \emph{partial order} in a set \(X\) is a relation `\(\sqsubseteq\)' between elements of \(X\) that is:
\begin{itemize}
    \item Extensional: \(x \sqsubseteq y \leftrightarrow (\forall x' \equiv x,\ y' \equiv y) (x' \sqsubseteq y').\) In fact, we shall only consider relations that are extensional. The interpretation is that they depend only on the object that the variable refers to;
    \item Decidable: \((x \sqsubseteq y) \lor \lnot (x \sqsubseteq y)\) (not to confuse with \emph{linearity}). In fact, we shall only consider relations that are decidable;
    \item Reflexive: \(x \sqsubseteq x\);
    \item Transitive: \(x \sqsubseteq y \sqsubseteq z \to x \sqsubseteq z\);
    \item Antisymmetric: \(x \sqsubseteq y \sqsubseteq x \to x \equiv y\).
\end{itemize}
From now on, all the spaces that we consider are equipped with a partial order `\(\sqsubseteq\)', the \emph{specialization order}.\footnote{It is sometimes advisable for the specialization order to be merely a pre-order, but we shall not worry about this here.} Furthermore, all functions \(f : X \to Y\) between spaces that we consider are \emph{covariant}:
\[ x \sqsubseteq x' \to f(x) \sqsubseteq f(x'). \]

Given a space \(X,\) a subspace \(A \subseteq X\) is given by an extensional membership predicate ``\(x \in A\)'' defined for \(x \in X.\) A subspace \(K\) is \emph{downward-closed} if:
\[ (\forall x \in K) (x \in K \land y \sqsubseteq x) \to (y \in K). \]
Every subspace \(A\) has a downward-closure:
\[ A\!\downarrow = \{x \in X\ |\ (\exists y \in A) (x \sqsubseteq y)\}. \]
A subspace \(F\) is \emph{upward-closed} if:
\[ (\forall x \in F) (x \in F \land y \sqsupseteq x) \to (y \in F). \]
Every subspace \(A\) has an upward-closure:
\[ A\!\uparrow = \{x \in X\ |\ (\exists y \in A) (x \sqsupseteq y)\}. \]
The complement of a downward-closed subspace is upward-closed, and vice-versa. We can construct the space
\[ \kappa(X) = \{K \subseteq X\ |\ K \text{ downward-closed}\}, \]
with specialization order \(K \sqsubseteq K' \leftrightarrow K \subseteq K'\), which is called the \emph{space of parts} of \(X.\) Hence, parts may shrink, but never grow, as we move to some accessible world. Note that \(\kappa(X)\) is a \emph{distributive lattice}: the union of two downward-closed subspaces is downward-closed, and the intersection of two downward-closed subspaces is downward-closed. In particular, the maximum element of \(\kappa(X)\) is \(X\) itself and its minimum element is \(\varnothing.\)

We have an embedding \(X \to \kappa(X)\) via \(x \mapsto x\!\downarrow.\) Sometimes we consider the embedding as an inclusion, in which case we write the unions of elements of \(\kappa(X)\) as joins `\(\sqcup\)', and the intersections as meets `\(\sqcap\)'. Note that each part \(K \in \kappa(X)\) can be represented as a join of elements of \(X\):
\[ K \equiv \bigsqcup_{x \in K} x\!\downarrow. \]
A \emph{prime} of \(\kappa(X)\) is an element \(x \in \kappa(X)\) such that \(x \not\equiv\varnothing\) and, for any other two elements \(K, K' \in \kappa(X)\):
\[ x \sqsubseteq K \sqcup K' \to (x \sqsubseteq K) \lor (x \sqsubseteq K'). \]
Equivalently, a prime is an element \(x \in \kappa(X)\) that satisfies, for any function \(I \to \kappa(X)\), \(i \mapsto K_i\):
\[ x \sqsubseteq \bigsqcup_{i \in I} K_i \to (\exists i \in I) x \sqsubseteq K_i. \]
The space of primes of \(\kappa(X)\) is precisely \(X\) itself. A function \(f : \kappa(X) \to \kappa(Y)\) is \emph{join-preserving} if:
\[ f\left(\bigsqcup_{i \in I} K_i \right) \equiv \bigsqcup_{i \in I} K_i. \]
Equivalently, \(f\) preserves joins if \(f(\varnothing) \equiv \varnothing\) and \(f(K \sqcup K') \equiv f(K) \sqcup f(K').\) It turns out that any function \(X \to \kappa(Y)\) (remember that it is assumed to be covariant) may be extended uniquely to a function \(\kappa(X) \to \kappa(Y)\) that preserves joins. Furthermore, any distributive lattice is of the form \(\kappa(X),\) where \(X\) is its space of primes (see \autoref{thm:birkhooff} later).

We also consider \(\mathfrak f(X),\) the lattice of \emph{closed subspaces} of \(X,\) which are the upward-closed subspaces of \(X,\) once again ordered by inclusion. If \(X\) is a space, we may form the \emph{dual space} \(X^* = \{x^*\ |\ x \in X\},\) where the specialization order is the reverse of that of \(X\):
\[ x^* \sqsubseteq y^* \leftrightarrow y \sqsubseteq x. \]
In particular, \(x^{**} \equiv x\) and \(X^{**} \equiv X.\)

A \emph{parent-child relation} between spaces \(X, Y\) is a part \(R \in \kappa(X^* \times Y).\) Equivalently, it is a relation \(R \subseteq X \times Y\) with the following property:
\[ x' \sqsupseteq x R y \sqsupseteq y' \to x' R y'. \]
If \(xRy,\) we say \(x\) is a \emph{parent} of \(y,\) or \(y\) is a \emph{child} of \(x.\) A parent-child relation may be understood as a way to extend the specialization order to elements of different spaces; it even makes sense to write `\(\sqsupseteq_R\)' for \(R\). Note that we have an \emph{identity} parent-child relation in each space, defined by \(xRx' \leftrightarrow x \sqsupseteq x'\). Furthermore, we can compose parent-child relations \(R : X \sqsupseteq Y \) and \(S : Y \sqsupseteq Z\) to obtain the parent-child relation \(SR : X \sqsupseteq Z\) given by:
\[ x SR z \iff (\exists y \in Y) (x R y S z). \]
Composition of parent-child relations is associative, and composition with the identity satisfies \(RI \equiv IR \equiv R\). Hence, partial orders with parent-child relations form a category. Furthermore, each parent-child relation \(R \in \kappa(X^* \times Y)\) comes with a dual relation \(R^* \in \kappa(Y \times X^*),\) defined by \(y^* R^* x^* \iff x R y.\)

Given a parent-child relation \(R \subseteq X \times Y\), we have a function \(R_* : \kappa(X) \to \kappa(Y)\), that maps each \(K \in \kappa(X)\) to its children:
\[ R_*(K) \equiv \bigsqcup_{\substack{x \in K\\y \in Y\\xRy}} y. \]
Reciprocally, we have a function \(R^{-1} : \mathfrak f(Y) \to \mathfrak f(X)\) that maps each \(F \in \mathfrak f(Y)\) to its parents:
\[ R^{-1}(F) \equiv \bigcup_{\substack{x \in X\\y\in F\\xRy}}x\!\uparrow. \]
In fact, \(R^{-1}\) is dual to \(R_*\), in the following sense:
\[ R^{-1}(F) \equiv ((R^*)_*(F^*))^*. \]

\begin{theorem}
    The following categories (whose objects are finite) are equivalent:
    \begin{itemize}
        \item The category of partial orders with parent-child relations;
        \item The category of join-preserving functions between distributive lattices.
    \end{itemize}
\end{theorem}

A \emph{parent-child relation} \(R \subseteq X \times Y\) is \emph{functional} if, for each \(x \in X\), there exists \(y_x \in Y\) such that, for all \(y \in Y\):
\[ x R y \leftrightarrow y_x \sqsupseteq y. \]
A functional parent-child relation gives rise to a function \(f_R : X \to Y,\) defined precisely by \(f_R(x) \sqsupseteq y \iff x R y.\) Reciprocally, given a function \(f : X \to Y\), we have a relation \(R_f \subseteq X \times Y\) defined by:
\[ x R_f y \iff f(x) \sqsupseteq y, \]
which is a parent-child relation because \(f\) is covariant. These constructions are inverses of each other. Furthermore, in this case, \(R^{-1}\) coincides with \(f^{-1} : \mathfrak f(Y) \to \mathfrak f(X)\), and so in particular \(f^{-1}(Y) \equiv X\) and \(f^{-1}(F \cap F') \equiv f^{-1}(F) \cap f^{-1}(F')\). That is, \(f^{-1}\) preserves \emph{meets} as well as joins.

\begin{theorem}[Birkhoff's representation theorem; see {\cite[p. 290]{lattice}}]
\label{thm:birkhooff}
    The following categories (whose objects are finite) are equivalent:
    \begin{itemize}
        \item The category of partial orders with functional parent-child relations between them.
        \item The category of partial orders with monotonic functions between them;
        \item The category of classical (again, finite) topological spaces;
        \item Dually, the category of distributive lattices with functions that preserve meets and joins.
    \end{itemize}
\end{theorem}

An important relation between variables is the following: \(x, z\) are \emph{indistinguishable}  if they have a lower bound in the specialization order:
\[ x \sim z \leftrightarrow (\exists y) (y \sqsubseteq z) \land (y \sqsubseteq x). \]
If \(x, y\) are not indistinguishable, they are said to be \emph{apart}, \(x \ne y\). We reserve `\(\ne\)' for apartness, rather than negated equality.

\section{The logic iS4}
\label{sec:logic}

The logic presented here diverges from the classical logic most often used in mathematical practice in two ways:
\begin{itemize}
    \item First, the system TN is based on \emph{intuitionistic logic}, which does not accept the principle of excluded middle \(\phi \lor \lnot \phi,\) or equivalently double-negation elimination \(\lnot\lnot\phi \to \phi,\) as valid. This does \emph{not} mean that one can take any particular instance of its negation as valid (in fact, \(\lnot (\phi \lor \lnot \phi) = \lnot \phi \land \lnot \lnot \phi\) is a contradiction just as much as it is in classical logic). However, a difference between TN and other constructive systems is that the most basic sentences in TN are \emph{judgements}, and hence decidable (that is, the law of excluded middle holds of judgements).
    \item Second, the system TN is based on \emph{modal logic}, which extends the language of logic using the two modalities `\(\Diamond\)' (\emph{possibly}) and `\(\square\)' (\emph{necessarily}). A sentence \(\phi\) that satisfies, necessarily, \(\phi \to \square \phi\) is said to be \emph{open} and a sentence that satisfies, necessarily, \(\Diamond \phi \to \phi\) is said to be \emph{closed}. The negation of a closed sentence is open and vice-versa, but again, the double negation \(\lnot\lnot\phi\) is not equivalent to \(\phi\) unless \(\phi\) is decidable. We shall talk more about these two symbols after introducing the axioms.
\end{itemize}
Formally, adding quantifiers to a modal theory is often not a trivial matter. Here we don't attempt to solve this issue: rather, we approach quantifiers informally, attempting to clearly describe the most important properties that the quantifiers satisfy. The modalities have proved not to be indispensable for developing the theory of spaces, if types are taken to be actually infinite. Nevertheless, it is more in the spirit to consider types to be merely potentially infinite, too.

\begin{axiom}[IS4, first-order version]
The axioms and rules of first-order intuitionistic logic are valid in TN, alongside the following modal axioms (see \cite{10.1093/logcom/exae030} for a similar presentation of the propositional version):
\begin{itemize}
\item \(\phi \to \Diamond \phi\) and \(\Diamond \Diamond \phi \to \Diamond \phi\);
\item \(\square \phi \to \phi\) and \(\square \phi \to \square \square \phi\);
\item \(\square (\phi \to \psi) \to (\square \phi \to \square \psi)\);
\item \(\square (\phi \to \psi) \to (\Diamond \phi \to \Diamond \psi)\);
\item \(\lnot \Diamond (\phi \land \lnot \phi)\);
\item \(\Diamond (\phi \lor \psi) \to \Diamond \phi \lor \Diamond \psi\);
\item \((\Diamond \phi \to \square \psi) \to \square (\phi \to \psi);\)
\item \((\exists x) \Diamond \phi(x) \leftrightarrow \Diamond (\exists x) \phi(x)\);
\item \((\exists x)\square \phi(x) \to \square (\exists x) \phi(x)\);
\item \(\Diamond (\forall x) \phi(x) \to (\forall x) \Diamond \phi(x)\);
\item \((\forall x) \square \phi(x) \leftrightarrow \square (\forall x) \phi(x)\);
\item \((\forall x)(\phi(x) \lor \lnot \phi(x)) \to ((\forall x) \phi(x) \lor (\exists x) \lnot \phi(x))\) (\emph{principle of omniscience}).
\end{itemize}
The necessitation rule to deduce \(\square \phi\) if \(\phi\) has been deduced (with no intermediate hypotheses) is valid, too.
\end{axiom}

As an interpretation of the modalities in terms, we imagine a fixed \emph{Kripke model}. This is a pre-order \((W, \le)\) of \emph{accessibility} between \emph{possible worlds} \cite[p. 34]{modal}, together with a relation ``holds'' (or \emph{``satisfies''}) \(w \Vdash \phi\) between worlds \(w \in W\) and sentences \(\phi\) such that \(w_0 \Vdash \square \phi\) if and only if, necessarily, for all \(w_1 \ge w_0,\) \(w_1 \Vdash \phi,\) and \(w_0 \Vdash \Diamond \phi\) if and only if there possibly exists \(w_1 \ge w_0\) such that \(w_0 \Vdash \phi.\) When we are studying some world \(w_0,\) we say that something happens \emph{later} if it happens in some world \(w_1 \ge w_0.\) Other temporal analogies should also be interpreted in terms of possible worlds, even if the order between worlds is not linear. We may now interpret each variable to be a filter; that is, each variable becomes increasingly specific.

\section{Types}
\label{sec:types}

A \emph{type} is like a space of spaces, with a specialization order between them. However, there are multiple ways for a space to be \emph{more specific} than another.\footnote{The inspiration here is \emph{Homotopy type theory} \cite{hottbook}, in which there are multiple ways that two elements can be equivalent.} Each of these is going to be a parent-child relation \(R : X \sqsupseteq Y.\) Hence, first we consider the type \(\mathcal E\) of all spaces, which is the category of all existing spaces and all parent-child relations between them. A \emph{type} is then defined to be an arbitrary subcategory of \(\mathcal E.\)\footnote{We could have followed a similar route with spaces and define the \emph{universe} \(\mathcal U\) to be the space of all variables where all of them are considered equivalent. Then a \emph{space} is simply a subspace of \(\mathcal U,\) with a finer specialization order.} We prefer to regard \(\mathcal E\) as \emph{potentially} infinite, rather than \emph{actually} infinite, in order to leave the door open to defining ``higher topologies'', even if we don't do this here. For our purposes, types can only \emph{grow}, that is, every object or morphism in a type \emph{necessarily} is in that type.

We often use the following naming convention for types: the objects (spaces) of a type \(X\) are denoted by sub-indices \(X_i \in X\), and the variables in \(X_i\) are denoted with the same sub-index \(x_i \in X_i.\) The arrows (parent-child relations) \(R \in X\) are denoted by:
\[ R : X_0 \sqsupseteq X_1. \]
That is, as a space (a filter with respect to the specialization order of the type) becomes more specific, its subindex increases, so that a subindex may be understood as an accessible world. Given a type \(X\), we may construct the \emph{pointed type} \(X_{p},\) whose objects are the (dependent) pairs \((X_i, x_i)\) where \(x_i \in X_i \in X\), and whose morphisms \(R : (X_0, x_0) \sqsupseteq (X_1, x_1)\) are the parent-child relations \(R : X_0 \sqsupseteq X_1\) such that \(x_0 R x_1.\) The \emph{points} of \(X\) are the filters of \(X_p\).

Given a type \(X,\) we may form the \emph{dual type} \(X^*,\) whose objects are the dual spaces \(X_i^*\) for \(X_i \in X\) and whose morphisms are \(R^* : X_1^* \sqsupseteq X_0^*\) for \(R : X_0 \sqsupseteq X_1\) a morphism of \(X\).

The most important properties that a morphism \(R : X_0 \sqsupseteq X_1\) may have, apart from the already mentioned \emph{functionality}, are:
\begin{itemize}
    \item \(R\) is \emph{overt} if children always exist: for all \(x_0 \in X_0,\) there exists \(x_1 \in X_1\) such that \(x_0 R x_1;\) equivalently, if \(R_*\) maps inhabited parts to inhabited parts, or \(R^{-1}(X_1) \equiv X_0\);
    \item \(R\) is \emph{compact} if parents always exist: for all \(x_1 \in X_1,\) there exists \(x_0 \in X_0\) such that \(x_0 R x_1\); equivalently, if \(R_*(X_0) \equiv X_1,\) or \(R^{-1}\) maps inhabited closed subspaces to inhabited closed subspaces;
    \item \(R\) is \emph{discrete} if any two children \(x_0 R x_1\), \(x_0 R x_1'\) have a third child \(x_0 R \hat x_1\) as an upper bound:
    \[ x_1 \sqsubseteq \hat x_1 \sqsupseteq x_1'. \]
    \item \(R\) is \emph{Hausdorff} if any two parents \(x_0 R x_1\), \(x_0' R x_1\) have a third parent \(\hat x_0 R x_1\)  as a lower bound:
    \[ x_0 \sqsupseteq \hat x_0 \sqsubseteq x_0'. \]
\end{itemize}
The composition of overt (resp. compact, discrete, Hausdorff) is overt (resp. compact, discrete, Hausdorff). The identity of any space is overt, compact, discrete and Hausdorff.  The dual of an overt relation is compact and vice-versa. The dual of a discrete relation is Hausdorff, and vice-versa. From this, we define the most important properties that a type \(X\) can possess:
\begin{itemize}
    \item \(X\) is \emph{overt} if all of its morphisms are overt;
    \item \(X\) is \emph{compact} if all of its morphisms are compact;
    \item \(X\) is \emph{discrete} if all of its morphisms are discrete;
    \item \(X\) is \emph{Hausdorff} if all of its morphisms are Hausdorff;
    \item \(X\) is \emph{directed} if, for all \(X_i, X_j \in X,\) there exists \(X_k \in X\) with \(X_i \sqsupseteq X_k \sqsubseteq X_j\);
    \item \(X\) is \emph{concrete} if, for all \(X_i, X_j \in X,\) there exists at most one morphism \(X_i \sqsupseteq X_j\) (in which case, its objects are called \emph{instances}).
\end{itemize}
In particular, the dual of an overt type is compact and vice-versa. The dual of a discrete type is Hausdorff, and vice-versa. Note that the discrete and Hausdorff conditions are satisfied by types such as the Sierpinski space. Hence, we might want to require stronger separation axioms.

Given a type \(X,\) we define the type \(\kappa(X)\) of \emph{parts} of \(X\), whose objects are (dependent) pairs \((X_i, K_i)\) with \(K_i \in \kappa(X_i)\),\footnote{Here, and for similar constructions, \(K_i\) is regarded as the space; that is, \(x_i \in (X_i, K_i) \iff x_i \in K_i\). The \emph{ambient space} \(X_i\) is regarded as additional information that \(K_i\) carries.}, and whose morphisms \(R : (X_0, K_0) \sqsupseteq (X_1, K_1)\) are the parent-child relations \(R : X_0 \sqsupseteq X_1\) in \(X\) such that \(R_*(K_0) \supseteq K_1,\) or equivalently:
\[ (\forall x_1 \in K_1) (\exists x_0 \in K_0) (x_0 R x_1). \]
In particular, the type of parts of an \emph{arbitrary} type is compact. This is one of the reasons that we expect TN to model locally compact spaces. We define the type \(\mathfrak f(X)\) of \emph{closed subspaces} of \(X\) to consist of the (dependent) pairs \((X_i, F_i)\) with \(F_i \in \mathfrak f(X_i)\), and whose morphisms \(R : (X_0, F_0) \sqsupseteq (X_1, F_1)\) are the parent-child relations \(R : (X_0, U_0)\) such that \(F_0 \sqsupseteq R^{-1}(F_0),\) or equivalently:
\[ (\forall x_1 \in F_1) (\forall x_0 \in X_0) (x_0 R x_1 \to x_0 \in F_0). \]
We may also define the \emph{open subspaces} to be complementary to the closed subspaces. Hence, \(\tau(X)\) is constructed as the type that consists of the pairs \((X_i, U_i)\) for \(U_i \sqsubseteq X_i \in X\) a downward-closed subspace (but now called an \emph{open subspace} rather than a part), and whose morphisms \(R : (X_0, U_0) \sqsupseteq (X_1, U_1)\) are the parent-child relations \(R : X_0 \sqsupseteq X_1\) such that \(R_*(U_0) \subseteq U_1,\) or equivalently:
\[ (\forall x_0 \in U_0) (\forall x_1 \in X_1) (x_0 R x_1 \to x_1 \in U_1). \]
A further characterization is:
\[ X_0 \setminus U_0 \supseteq R^*(X_1 \setminus U_1), \]
which is why we say that closed and open subspaces are complementary. In fact, the category \(\tau(X)\) is equivalent to \(\mathfrak f(X)\), but since we care about the elements too, the correct equivalence of types involves two dualities:
\[ \tau(X) \cong \mathfrak f(X^*)^*. \]

We already said that parts are always compact. For open and closed subspaces, the situation is only slightly more complex:
\begin{theorem}
    If \(X\) is overt, then \(\tau(X)\) is overt. If \(X\) is compact, then \(\mathfrak f(X)\) is compact.
\end{theorem}
\begin{proof}
    Suppose \(X\) is an overt type, and let \(R : (X_0, U_0) \sqsupseteq (X_1, U_1)\) be a morphism. Since \(R\) is overt, for all \(x_0 \in U_0\) there exists \(x_1 \in X_1\) such that \(x_0 R x_1.\) Since \(U\) is open, this implies \(x_1 \in U_1.\)

    For the dual, we could carry out basically the same proof, but there is an even simpler way: if \(X\) is compact, then \(X^*\) is overt, so that \(\tau(X^*)\) is overt, so that \(\tau(X^*)^* \cong \mathfrak f(X)\) is compact.
\end{proof}

More generally, one can form the type \(\mathcal P(X)\) of all subspaces, which consists of the dependent pairs \((X_i, A_i)\) for \(A_i \subseteq X_i \in X\), where \(A_i\) is the space and \(X_i\) is extra information, and whose morphisms \(R : (X_0, A_0) \sqsupseteq (X_1, A_1)\) are all the morphisms \(R : X_0 \sqsupseteq X_1\) in \(X.\) A \emph{type of subspaces} is a subtype of \(\mathcal P(X),\) obtained by choosing some subspaces and some morphisms. Given a type of subspaces \(\mathcal A,\) one can form the \emph{downward-closure} \(\mathcal A\!\downarrow\), consisting of the \((X, A\!\downarrow)\) for \(A \in \mathcal A\), with precisely the same morphisms as \(\mathcal A\), and similarly the upward-closure \(\mathcal A\!\uparrow\). We say that a subspace \((X, A) \in \mathcal A\) \emph{has all upper bounds} if, for all \(a, a' \in A\), \(x \in X,\) we have:
\[ a \sqsubseteq x \sqsupseteq a' \implies x \in A. \]
If all subspaces have this property, then \(\mathcal A\) is said to have all upper bounds. Similarly, we can define types with all \emph{lower} bounds.

\begin{theorem}
    If \(X\) is discrete and \(\mathcal A \subseteq \mathcal P(X)\) has all upper bounds, then \(\mathcal A\) is discrete. If \(\mathcal A\) is also downward-closed and overt, then \(\mathcal A\) is \emph{open}:
    \[ \mathcal A \subseteq \tau(X). \]
    If \(X\) is Hausdorff and \(\mathcal A\) has all lower bounds, then \(\mathcal A\) is Hausdorff. If \(\mathcal A\) is also upward-closed and compact, then \(\mathcal A\) is closed:
    \[ \mathcal A \subseteq \mathfrak f(X). \]
\end{theorem}
\begin{proof}
    Suppose \(X\) is discrete and \(\mathcal A\) and has all upper bounds. Let \(R : (X_0, A_0) \sqsupseteq (X_1, A_1)\) be a morphism in \(\mathcal A,\) and let \(a_0 \in A_0\) have two children in \(A_1.\) Since \(X\) is discrete, the two children have an upper bound \(\hat a_1 \in X_1\). Since \(A\) has all upper bounds, \(\hat a_1 \in A_1.\) This shows that \(R\) is discrete.

    Suppose that, furthermore, \(\mathcal A\) is downward-closed and overt. Given \(a_0 \in A_0,\) and a child \(a_0 R x_1 \in X_1,\) we want to show that \(x_1 \in A_1.\) Since \(A\) is overt, there is some other child \(a_0 R a_1 \in A_1.\) Since \(A\) is discrete, these two children have an upper bound in \(A_1.\) Hence, \(x_1\) is above an element of \(A_1,\) but since \(A_1\) is downward-closed, this shows that \(x_1 \in A_1.\)

    Finally, if \(X\) is Hausdorff and \(\mathcal A\) has all lower bounds, then \(X^*\) is discrete and \(\mathcal A^*\) has all upper bounds, so that \(\mathcal A^*\) is discrete, that is, \(\mathcal A\) is Hausdorff. If \(\mathcal A\) is also upward-closed and compact, then \(\mathcal A^*\) is downward-closed and overt, so that \(\mathcal A^*\) is open; that is, \(\mathcal A\) is closed.
\end{proof}

\section{The Sierpinski space}
\label{sec:sierpinski}

Let \(\Sigma = \{\mathbb O \sqsubseteq \mathbb F\}\) be the \emph{Sierpinski space}. See, for instance, \cite{TaylorP:lamcra} for an explanation. Three alternative descriptions of \(\Sigma\) are the following:
\begin{itemize}
    \item \(\Sigma \equiv \kappa(\{*\}) = \{\varnothing \sqsubseteq \{*\}\}\) is the lattice of parts of a single point;
    \item \(\Sigma \equiv \mathfrak f(\{*\}) = \{\varnothing \sqsubseteq \{*\}\}\) is the lattice of closed subspaces of a single point;
    \item \(\Sigma \equiv \tau(\{*\}) = \{\{*\} \sqsubseteq \varnothing\}\) is the lattice of open subspaces of a single point.
\end{itemize}
\(\mathbb O\) is the \emph{open point} (since \(\{\mathbb O\}\) is the only non-trivial open subspace) and \(\mathbb F\) is the \emph{closed point}.
\begin{theorem}
    Given a type \(X,\) the following 4 categories are equivalent:
    \begin{itemize}
        \item The category of parent-child relations \(\chi_i : X_i \sqsupseteq \{*\}\) for \(X_i \in X\), called \emph{characteristic maps};
        \item The category of functions \(f_i : X_i \to \Sigma\) for \(X_i \in X\), called \emph{classifying maps};
        \item The category \(\mathfrak f(X)\) of \emph{closed subspaces}.
        \item The category \(\tau(X)\) of \emph{open subspaces}.
    \end{itemize}
\end{theorem}
\begin{proof}
    We need to define the first two categories. Given characteristic maps \(\chi_0 : X_0 \sqsupseteq \{*\},\) \(\chi_1 : X_1 \sqsupseteq \{*\}\), the morphisms between \(\chi_0\) and \(\chi_1\) are defined to be the parent-child relations \(R : X_0 \sqsupseteq X_1\) such that \(\chi_0 \sqsupseteq \chi_1 R.\) A morphism between classifying functions \(R : (X_0, f_0) \sqsupseteq (X_1, f_1)\) is a morphism \(R : X_0 \sqsupseteq X_1\) such that, for all \(x_0 R x_1\), we have \((f_0(x_0) \sqsubseteq \mathbb O) \to (f_1(x_1) \equiv \mathbb O),\) or equivalently, \(f_1(x_1) \equiv \mathbb F \to f_0(x_0) \equiv \mathbb F.\) (This is a particular case of \emph{covariance} between multiple spaces). We also recall that \(\tau(X)\) is equivalent to \(\mathfrak f(X)\) (as \emph{categories}, not as \emph{types}).
    
    Given a characteristic map \(\chi_i\), \(\chi_i^{-1}(\{*\}) = \{x \in X\ |\ x \chi *\}\) is upward-closed (called the \emph{classified closed subspace}). Let \(R : (X_0, \chi_0) \sqsupseteq (X_1, \chi_1)\) be a morphism of characteristic maps. If \(x_0 \in X_0,\ x_1 \in \chi_1^{-1}(\{*\})\) are such that \(x_0 R x_1\), then \(x_0 R x_1 \chi_1 *\). Since \(\chi_0 \sqsupseteq \chi_1 R),\) we have \(x_0 \chi_0 *,\) or \(x_0 \in \chi_0^{-1}(\{*\}).\) This shows that
    \[ R : (X_0, \chi_0^{-1}(\{*\})) \sqsupseteq (X_1, \chi_0^{-1}(\{*\})) \]
    is a morphism of closed subspaces.

    Conversely, given a closed subspace \((X_i, F_i)\), define the characteristic map \(x_i \chi_i * \iff x_i \in F_i.\) It is a parent-child relation because \(F_i\) is upward-closed; furthermore, it is the inverse of the above construction. Given a morphism \(R : (X_0, F_0) \sqsupseteq (X_1, F_1)\) of closed subspaces:
    \[ x_0 \chi_0 * \iff x_0 \in F_0 \impliedby (\exists x_1 \in F_1) (x_0 R x_1) \iff (\exists x_1)\, x_0 R x_1 \chi_1 *. \]
    This shows that \(\chi_0 \sqsupseteq \chi_1 R.\)

    Given a characteristic map \((X_i, \chi_i)\), we have the function \( \chi_{i*} : \kappa(X_i) \to \kappa(\{*\}) \equiv \Sigma \), whose restriction to \(X_i \subseteq \kappa(X_i)\) is the desired classifying function. If \(R : (X_0, \chi_0) \sqsupseteq (X_1, \chi_1)\) is a morphism of classifying maps, and \(\chi_*(x_1) \equiv \mathbb F\) for some \(x_1\) (that is, \(x_1 \chi_1 *\)), then for any \(x_0 R x_1\) we have \(x_0 R x_1 \chi_1 *\), so that, in particular, \(x_0 \chi_0 *.\)
    
    Finally, a characteristic map is obtained from a classifying function \(f_i\) by \(x_i \chi_i * \iff f_i(x_i) \equiv \mathbb F.\) Again, the construction is inverse to the one above. If \(R : (X_0, f_0) \sqsupseteq (X_1, f_1)\) is a morphism of classifying functions, and \(x_0 R x_1 \chi_1 *,\) then \(f_1(x_1) \equiv \mathbb F,\) so that \(f_0(x_0) \equiv \mathbb F;\) that is, \(x_0 \chi_0 *,\) which shows that \(R : (X_0, \chi_0) \sqsupseteq (X_1, \chi_1)\) is a morphism of characteristic maps.
\end{proof}

\section{Functions between types}
\label{sec:functions}

For specific functions, the domain and codomain are usually \emph{directed, concrete types}. Hence, there is always a most specific instance, and each object specializes to it in a unique way. A \emph{function of type} \(f : X \to Y\) between concrete, directed types is given by:
\begin{itemize}
    \item A full subtype \(T\) of \(X \times Y\);\footnote{A full subtype contains all the morphisms between its objects. In this case, this just means that the specialization order of \(T\) is the restriction of the specialization order of \(X \times Y\).}
    \item For each \((X_i, Y_j) \in T\), a function \(f_{ij} : X_i \to Y_j\);
    \item For each \((X_0, Y_0) \in X \times Y,\) there (possibly) exists \((X_0, Y_0) \sqsupseteq (X_1, Y_1) \in T\);
    \item \emph{Filter covariance}: If \(R^X : X_0 \sqsupseteq X_1\), \(R^Y : Y_0 \sqsupseteq Y_1,\) and \((X_0, Y_0),\ (X_1, Y_1) \in T\), then, whenever \(x_0 R^X x_1\):
    \[ f_{00}(x_0) R^Y f_{11}(x_1). \]
\end{itemize}

\section{Rational numbers}
\label{sec:rationals}

The type of unbounded dense linear orders consists of the discrete spaces \(Q\) with a linear order `\(<\)'. Its morphisms are the order-preserving injective functions; in particular, they are overt. As a consequence of the definition, a generic dense linear order \(Q\) (i.e., a filter of this type, with no additional assumptions) satisfies:
\begin{equation}
\label{eq:dense}
a < c \in Q \implies \Diamond (\exists b \in Q) (a < b < c).
\end{equation}
Also by definition of an \emph{unbounded} order:
\begin{equation}
\label{eq:unbounded}
(\forall b \in Q) \Diamond (\exists a, c \in Q) (a < b < c).
\end{equation}
(The type of bounded dense linear orders has the same spaces, but its morphisms are the order-preserving injective functions that preserve maxima and minima.)

Rather than constructing the rationals, we suppose we are given a concrete, directed subtype \(\mathbb Q\) with the following properties:
\begin{itemize}
\item The specialization order is given by reverse inclusion;
\item It is still unbounded and dense in the sense of \eqref{eq:unbounded} and \eqref{eq:dense};
\item Every instance shares the point \(0;\) 
\item An addition function of type \(\mathbb Q \times \mathbb Q \to \mathbb Q\), together with a negation operation of type \(\mathbb Q \to \mathbb Q\), with the usual properties. For example, associativity is stated as \((x+y)+z \equiv x+(y+z)\) for any choice of instances where the required additions are defined. Directedness ensures that the equation always make in a specific enough (large enough) instance, and concreteness ensures that it makes sense in a single way, provided the left-hand side and right-hand side are defined somewhere. In this case, once addition is defined, the domain may be replaced by any more general domain, and the codomain by any more specific codomain.
\end{itemize}

\section{Ascending and descending line}
\label{sec:ascending}

The type of \emph{ascending lines} (see \cite{TaylorP:dedras}) consists of the spaces \(R_{\nearrow}\) that carry a linear order `\(\le\)' that is opposite to their inclusion order:
\[ \delta \sqsubseteq \delta'  \iff \delta \ge \delta'. \]
The morphisms between ascending lines are the strictly increasing functions between them. (Note that this is the same as covariance.) Again, the ascending line is overt.

Symmetrically, the descending line consists of the spaces with a linear order that coincides with the specialization order:
\[ \upsilon \sqsubseteq \upsilon' \iff \upsilon \le \upsilon', \]
with (covariant) functions as morphisms.

We also have bounded versions of these. The ascending line type with a lower bound is compact, since each ascending element is a child of the minimum, and so is the descending line with an upper bound.

As particular subtypes, we have \(\mathbb R_{\nearrow}\) and \(\mathbb R_{\swarrow}\), which are directed types with one instance \(\mathbb R^i_{\nearrow}\) and one instance \(\mathbb R^i_{\swarrow}\) for each instance \(\mathbb Q^i,\) and whose elements are \(a_{\nearrow} \in \mathbb R_{\nearrow}\) and \(a_{\swarrow} \in \mathbb R_{\swarrow}\) for \(a \in \mathbb Q,\) with the same comparison order.

However, the morphisms of \(\mathbb R_{\nearrow}\) are the strictly increasing functions \(f: \mathbb R^0_{\nearrow} \to \mathbb R^1_{\nearrow}\) such that \(f(x) \ge x\) (i.e., \(f(x) \sqsubseteq x\)) rather than the stronger \(f(x) \equiv x\) that rationals satisfy. Similarly, the morphisms of \(\mathbb R_{\swarrow}\) are the strictly increasing functions that satisfy \(f(x) \le x.\)

\section{Line}
\label{sec:line}

The type of \emph{lines} consists of the spaces of the form \(R \subseteq R_{\nearrow} \times R_{\swarrow}\), where \(R_{\nearrow}\) and \(R_{\swarrow}\) are a ascending and a disjoint descending line, with the following properties:\footnote{The construction is similar to the one in \cite{TaylorP:dedras}. It is a curious fact that the line thus define boils down to a graph with increasingly many edges. Disjointness is just a stylistic preference; it will not affect the directedness of \(\mathbb R\).}
\begin{itemize}
\item \(R_{\nearrow} \cup R_{\swarrow}\) is still linearly ordered, with minimum \(\min R = \min R_{\nearrow}\) and maximum \(\max R = \max R_{\swarrow}\), and furthermore (for convenience) the elements of \(R_{\nearrow}\) always alternate between an ascending and a descending element. The picture is thus like this:
\[ R_{\nearrow} \cup R_{\swarrow} = \{\delta^0 \prec \upsilon^0 \prec \delta^1 \prec \upsilon^1 \prec \cdots \prec \delta^s \prec \upsilon^{s} \}, \]
where \(\delta^i \in R_{\nearrow}\) and \(\upsilon^i \in R_{\swarrow}.\)
\item Its elements are the \emph{primes}:
\[ (\delta, \upsilon) \in R \iff (\delta < \upsilon) \land \lnot(\exists \delta' \in R_{\nearrow}, \upsilon' \in R_{\swarrow}) (\delta < \delta' < \upsilon' < \upsilon).\]
This is to distinguish the line from the \emph{interval domain}, or so that joins work as one would expect;
\item The specialization order is that of the product, that is:
\[ (\delta, \upsilon) \sqsubseteq (\delta', \upsilon') \iff \delta' \le \delta < \upsilon \le \upsilon'; \]
\item The morphisms are given by order-preserving functions \(f : R_{\nearrow}^0 \cup R_{\swarrow}^0 \to R_{\nearrow}^1 \cup R_{\swarrow}^1\), that restrict to a morphism of ascending elements and a morphism of descending elements.
\end{itemize}

The \emph{vertices} of the line are the points of the form \((\delta, \upsilon)\) for \(\delta\) ascending and \(\upsilon\) descending with \(\delta \prec \upsilon\). (This notation means that \(\delta < \upsilon\) and there is no element between them). Note that vertices shrink over time. When the ascending-descending pair is ordered backwards, \(\upsilon \prec \delta\), then \((\upsilon, \delta)\) is called an \emph{edge}, but an edge is not a point of the line, and it \emph{grows} rather than shrink. The points of the line that are not vertices are the \emph{closed edges} \((\delta, \upsilon)\) for \(\delta \prec \upsilon' \prec \delta' \prec \upsilon\); they are called this way because they are the union of the vertices \((\delta, \upsilon')\) and \(\delta', \upsilon\) with the edge \((\upsilon', \delta').\) A closed edge also shrinks over time.

A version of the line that is bounded below is obtained by using a version of the ascending numbers that is bounded below. Similarly, a version of the line that is bounded above is obtained by using a version of the descending numbers that is bounded above. The type of the bounded line (both above and below) is compact. Bounded or not, the line is overt and Hausdorff.

Given \(\delta \in R_{\nearrow}\), \(\upsilon \in R_{\swarrow}\) with \(\delta < \upsilon\), we may construct the closed interval \([\delta, \upsilon] = \{ (\alpha, \beta) \in R\ |\ (\alpha < \upsilon) \land (\delta < \beta) \},\) which is a closed subspace, and in fact:
\begin{theorem}
    Any closed interval \([\delta, \upsilon]\) is compact. If \(\delta, \upsilon \in R,\) then \([\delta, \upsilon]\) is overt, too.
\end{theorem}
\begin{proof}
    Let \(f : R_0 \to R_1\) be a morphism of lines. Since \(\delta\) is ascending, we have \(\delta_1 \ge f(\delta_0)\), and since \(\upsilon\) is descending, we have \(\upsilon_1 \le f(\upsilon_0).\) Suppose \((\alpha_1, \beta_1) \in [\delta_1, \upsilon_1],\) that is, \(\alpha_1 < \upsilon_1\) and \(\delta_1 < \beta_1.\) We need to find a parent for \((\alpha_1, \beta_1).\)
    
    Let \(\alpha_0\) be the largest ascending element in \(R_{\nearrow0}\) such that \(f(\alpha_0) \le \alpha_1,\) and let \(\beta_0\) be the smallest descending element in \(R_{\swarrow0}\) such that \(f(\beta_0) \ge \beta_1.\) If \((\alpha_0, \beta_0)\) were not prime, then \((\alpha_1, \beta_1)\) would not be prime either. Hence \((\alpha_0, \beta_0) \in R\) has \((\alpha_1, \beta_1)\) as a child. Since \(\alpha\) is ascending and \(\upsilon\) is descending:
    \[ \alpha_1 < \upsilon_1 \implies \alpha_0 < \upsilon_0, \]
    and since \(\delta\) is ascending and \(\beta\) is descending:
    \[ \delta_1 < \beta_1 \implies \delta_0 < \beta_0. \]
    In conclusion, \((\alpha_0, \beta_0) \in [\delta_0, \upsilon_0],\) which establishes compactness.

    Now suppose \(\delta \in R.\) What we mean by this is that \(\delta_1 \equiv f(\delta_0)\), so that we may assign the smallest element above \(\delta\) as a companion to make a vertex. In particular, any real element above \(\delta_0\) will remain above \(\delta_1.\) Suppose also \(\upsilon \in R,\) that is \(\upsilon_1 \equiv f(\upsilon_0),\) and take any \((\alpha_0, \beta_0) \in [\delta_0, \upsilon_0].\) We want to find a child \((\alpha_1, \beta_1)\) that remains in the closed interval.

    Since \(\beta_0 > \delta_0,\) we have that \(f(\beta_0) > \delta_1,\) by the discussion above. Similarly, since \(\alpha_0 < \upsilon_0,\) we have that \(f(\alpha_0) < \upsilon_1.\) Hence \((f(\alpha_0), f(\beta_0))\) is our candidate; the only problem is that it need not be prime. In any case, it \emph{contains} a prime. For example, define \(\alpha_1 \equiv \max\{f(\alpha_0), \delta_1\}\) and let \(\beta_1\) be the smallest element above \(\alpha_1.\) Then \(\alpha_1 \ge f(\alpha_0)\), and since \(f(\beta_0) > \alpha_1,\) it follows that \(\beta_1 \le f(\beta_0)\). Hence \(\alpha\) ascends and \(\beta\) descends, and furthermore \((\alpha_1, \beta_1)\) is a vertex, which satisfies \(\alpha_1 \le \max\{f(\alpha_0), \delta_1\} < \upsilon_1\) and \(\beta_1 > \alpha_1 \ge \delta_1,\) that is, \((\alpha_1, \beta_1) \in [\delta_1, \upsilon_1].\)
\end{proof}

In particular, one has the real numbers \(\mathbb R,\) which are the directed, concrete subtype obtained as the line from \((\mathbb R_{\nearrow}, \mathbb R_{\swarrow}).\) We shall construct addition. For any choice of three instances \(\mathbb R^0_0, \mathbb R^1_0, \mathbb R^2_0\), we shall find three more specific instances \(\mathbb R^i_0 \sqsupseteq \mathbb R^i_1\) where the addition function can be defined. The easiest way to achieve this is by assigning four numbers to each instance, that measure how `specific' the instance is:
\begin{itemize}
\item The first two numbers are the bounds \(\min \mathbb R^i\) and \(\max \mathbb R^i\). These are variable numbers: for example, the lower bound descends as new vertices are added to the left, but ascends as vertices shrink.
\item The third number is the \emph{Lebesgue number} \(\epsilon\), which is the smallest width of the vertices of \(\mathbb R^i,\) or equivalently:
\[ \epsilon(\mathbb R^i) = \min_{(\delta, \upsilon) \in \mathbb R^i} (\upsilon - \delta). \]
The Lebesgue number is a positive descending number. If \(a, b \in \mathbb Q\) are such that \(b-a \le \epsilon(\mathbb R^i)\) and \(a \ge \min \mathbb R^i,\) \(b \le \max \mathbb R^i,\) then \((a, b)\) is contained in some point of \(\mathbb R^i.\)
\item The last number is the \emph{precision} \(\delta\), which is the \emph{largest} width of a closed edge:
\[ \delta(\mathbb R^i) = \max_{(\delta, \upsilon) \in \mathbb R^i} (\upsilon - \delta). \]
The precision is a variable number. For the bounded version of \(\mathbb R,\) the precision is a descending number.
\end{itemize}

Suppose given \(\mathbb R^i_0\), \(i = 1,2,3\). Then we can construct an instance \(\mathbb R^3_1\) more specific than \(\mathbb R^3_0\) satisfying:
\[ \min \mathbb R^3_1 \le \min \mathbb R^1_0 + \min \mathbb R^2_0,\]
and the reverse inequality for the maximums. Let \(\epsilon\) be the Lebesgue number of \(\mathbb R^3_1\). Construct the number \(\delta = \epsilon/2\). Then we can find an instance \(\mathbb R^1_1\) more specific than \(\mathbb R^1_0\) with \(\min \mathbb R^1_1 \ge \min \mathbb R^1_0\), \(\max \mathbb R^1_1 \le \max \mathbb R^1_0\) and precision less than or equal to \(\delta.\) We construct \(\hat{\mathbb R}^2_1\) in the same way. Then we can define:
\begin{align*}
\hat{\mathbb R}^1 \times \hat{\mathbb R}^2 &\to \hat{\mathbb R}^3 \\
(a, b) + (a', b') &= (c, d),
\end{align*}
where \((c, d)\) is the smallest interval such that \( c \le a + a' < b+b' \le d \), which is a prime because \((b+b') - (a+a') = (b-a) + (b' -a') \le \delta + \delta \equiv \epsilon\). It turns out that the proof of continuity is not that different from the traditional ones. We sum up the construction of the addition operation between real numbers as follows:
\begin{itemize}
\item Addition is defined uniquely whenever we have three instances \(\mathbb R^1, \mathbb R^2, \mathbb R^3\) of the real numbers that satisfy \(\min \mathbb R^3 \le \min \mathbb R^2 + \min \mathbb R^3,\) \(\max \mathbb R^3 \ge \max \mathbb R^2 + \max \mathbb R^1\) and \(\delta(\mathbb R^1) + \delta(\mathbb R^2) \le \epsilon(\mathbb R^3)\).
\item For any three instances \(\mathbb R^1, \mathbb R^2, \mathbb R^3\) of the real numbers, we can construct more specific instances \(\hat{\mathbb R}^1, \hat{\mathbb R}^2, \hat{\mathbb R}^3\) that satisfy the three conditions above. This is done by making the codomain long enough (\emph{extending} the line) and by improving the precision of the domain enough, by shrinking its points (\emph{refining} the line).
\end{itemize}

\section{Applications}
\label{sec:applications}

In order to get a feel of how reasoning in TN and related areas (synthetic topology, ASD) works, this section provides a high-level overview of how to encode stuff from ordinary mathematics. Since the details have not been worked out at the moment, the aim is just to provide the intuitions that motivate TN.

First, let's review some places where ascending numbers appear in ordinary mathematics. While these are probably known in the folklore of synthetic topology, I believe they deserve attention by the wider mathematical audience:
\begin{itemize}
\item The rank of an \(n \times n\) matrix is an ascending natural number between \(0\) and \(n.\) This is because the determinant of a real matrix being non-zero is an open condition. It is sometimes said that the rank is not well-defined in practice, but it is more accurate (at least when doing exact real arithmetic rather than floating-point arithmetic) to say that the rank is ascending. A \(2\times2\) matrix whose rank is known to be above \(1\) is equivalent to a matrix of the form
\[\begin{pmatrix}
  1 & 0 \\ 0 & x
  \end{pmatrix},\]
where \(x\) is currently indistinguishable from \(0,\) but may later be found to actually be non-zero. However, note that the rank of a matrix with coefficients in a discrete field such as \(\mathbb Q\) is a natural number.
\item The supremum of an open subspace of the real line, and the infimum of a closed part of the real line, are ascending numbers. This is because an open subspace may grow and a closed part may shrink.
\item Given a point \(p\) and a tangent vector \(v\) in a Riemannian manifold, the domain of definition of the geodesic curve that crosses \(p\) with velocity \(v\) is an open interval. The domain may grow as the point shrinks and the direction becomes more specific, because a finer geodesic may no longer hit the `obstacle' that prevented its extension. Hence, the maximum length of the geodesic is an ascending real number.
\item The Lebesgue measure of an open subspace of \(\mathbb R^n\) is an ascending number.
\item The largest natural number (or the largest number that has been defined) is ascending.
\item The ``current world'' \(w_c\) where the theory is developed ascends, from the point of view of the theory itself (even if \(w_c\) can only be partially defined), and \(w_c\) \emph{usually}, but not always, ascends from the point of view of the meta-theory, too.
\end{itemize}

Now let's see how basic homotopy theory could work in TN. We can construct the circle \(S\) in a way analogous to the line, as a graph where each edge may be subdivided into two. This time the graph forms a loop rather than a line. There are at least 3 vertices in \(S\). For any instance of the circle, one can find an instance of the line as long as one needs, with a covering map \(\pi : R \to S.\) Then we would show:
\begin{theorem}
A continuous map \(\gamma : R \to S\), where \(R\) is a line, can be lifted to a continuous map \(\hat \gamma : R \to R'\), where \(\pi : R' \to S\) is a covering map and \(\pi \hat \gamma = \gamma.\) One endpoint of \(\hat \gamma\) determines the other, even if \(\gamma\) undergoes a homotopy with fixed endpoints.
\end{theorem}

It should also not be surprising if basic analysis could be done in TN, considering its relation to ASD. What is a differentiable (or \(\mathcal{C}^1\)) map \(f: \mathbb R \to \mathbb R\)? It should simply be a continuous map for which there is another continuous map of type \(F : \mathbb R \times \mathbb R \to \mathbb R\) such that \(F|_{(\mathbb R \times \mathbb R) \setminus \Delta}(x, y) = ((f(y) - f(x))/(y-x)\), where \(\Delta = \{(x, y) \in \mathbb R \times \mathbb R\ |\ x \sim y\}\) is the (closed) diagonal.

A continuous function \(f : \mathbb R^1 \to \mathbb R^2\) could also be integrated:
\[ \int f = \sum_{\substack{a \prec b \in \mathbb Q^1}} (b-a) f(a, b), \]
and a weight could also be added to the vertices, if necessary. A shortcoming is that the endpoints here are rational numbers (or members of some other \emph{discrete}, dense linear order). Apart from that, for integration purposes, it appears to make more sense to use a finer topology \emph{in the domain}: define \(\kappa(\mathcal L(R))\) to be the boolean algebra generated by \(R\), or to be more precise, \(\mathcal L(R)\) is the discrete topological space on the set of vertices \(a \in Q\) and edges \((a, b)\) with \(a \prec b\).

Whether the definitions in this section are viable in practice remains to be seen. In any case, the biggest obstacle TN might face is its extension to spaces that are not locally compact.

\textbf{Acknowledgements.} The author would like to thank Nik Weaver for a short exchange involving an early summary, in which he advised me to make the exposition more accessible.

\textbf{Artificial intelligence acknowledgement.} Artificial intelligence was not used to help redact this paper or to generate any of its ideas.


\providecommand{\bysame}{\leavevmode\hbox to3em{\hrulefill}\thinspace}
\providecommand{\MR}{\relax\ifhmode\unskip\space\fi MR }
\providecommand{\MRhref}[2]{%
  \href{http://www.ams.org/mathscinet-getitem?mr=#1}{#2}
}
\providecommand{\href}[2]{#2}

\end{document}